\documentclass{amsart}
\usepackage[croatian]{babel}

\usepackage{amsmath, amssymb, amsfonts}
\usepackage[all]{xy}
\usepackage{enumitem}   
\usepackage{float}

\newtheorem{tm}{Theorem}

\newenvironment{dokaz}
{\noindent\emph{Proof:}\ }
{\hfill $\blacksquare$}

\newcommand{\Z}
{{\mathbb Z}}
\newcommand{\N}
{{\mathbb N}}
\newcommand{\C}
{{\mathbb C}}
\newcommand{\m}
{{\mathbf m}}

\newcommand{\g}
{{\mathfrak g}}

\newcommand{\gt}
{\tilde{{\mathfrak g}}}

\newcommand{\h}
{{\mathfrak h}}
\newcommand{\n}
{{\mathfrak n}}
\newcommand{\gsl}
{{\mathfrak sl}}

\newcommand{\Gamt}
{\tilde{{\Gamma}}}

\begin{document}

\title[Presentations of Feigin-Stoyanovsky's type subspaces]{Presentations of Feigin-Stoyanovsky's type subspaces of standard modules for affine Lie algebras of type $C_\ell^{(1)}$}
%\Country{University of Zagreb, Croatia}
\author{Goran Trup\v{c}evi\'{c}}
\address{Faculty of Teacher Education, University of Zagreb, Zagreb, Croatia}
%%\curraddr{}
\email{goran.trupcevic@ufzg.hr}
\thanks{Partially supported by the Croatian Science Foundation grant 2634.}
%%\date{}
%%\dedicatory{}
%
\begin{abstract}
	Feigin-Stoyanovsky's type subspace $W(\Lambda)$ of a standard $\gt$-module $L(\Lambda)$ is a $\gt_1$-submodule of $L(\Lambda)$ generated by
	 the highest-weight vector $v_\Lambda$, 	where $\gt_1$ is a certain commutative subalgebra of $\gt$. Based on the description of basis of 
	 $W(\Lambda)$ for $\gt$ of type $C_\ell^{(1)}$, we 	give a presentation of this subspace in terms of generators and relations
	$$W(\Lambda)\simeq U(\gt_1^-)/J.$$
\end{abstract}

\maketitle

\section{Introduction}

B. Feigin and A. Stoyanovsky introduced principal subspaces of standard modules for affine Lie algebras of type $A_1^{(1)}$ and $A_2^{(1)}$ 
in \cite{FS} where they have recovered Rogers-Ramanujan type identities by considering graded dimensions of these subspaces. An important part
of their investigation was the knowledge of presentations of these subspaces in terms of generators and relations. 
Another type of principal subspaces, called  Feigin-Stoyanovsky's type subspace, was introduced by M. Primc who constructed bases of these subspaces
in different cases (\cite{P1}, \cite{P2}, \cite{P3}, \cite{JP}). 
These kind of subspaces were further studied by many authors (\cite{Sto}, \cite{G}, \cite{CLM1}, \cite{CLM2}, \cite{AKS}, \cite{FJLMM}, \cite{FJMMT}, \cite{Ba}, \cite{Bu1}, \cite{Bu2}, \cite{J1}, \cite{T1}, \dots) and the knowledge of presentation
presents an important question in this study (\cite{C1},  \cite{CalLM1}, \cite{CalLM2}, \cite{CalLM3}, \cite{S1}, \cite{Pe}).

In our previous works  we have described bases of Feigin-Stoyanovsky's type subspaces  of standard modules for affine Lie algebras 
of type $C_\ell^{(1)}$ (\cite{BPT}) and obtained from them basis for the whole standard modules (\cite{BPT}). In this note we use the description of bases of a Feigin-Stoyanovsky's type subspaces to give presentations of these subspaces in terms of generators and relations.

\section{Feigin-Stoyanovsky's type subspaces}

%%%%%%%%%%%%%%%%%%%%%%%%%%%%%%%%%%%%%%%%%%%%%%%%%%% prosta kondim Liejeva algebra

Let ${\mathfrak g}$ be a complex simple  Lie algebra of type $C_\ell$ with a Cartan
subalgebra $\h$ and a root decomposition ${\mathfrak g}={\mathfrak h}+\sum {\mathfrak g}_\alpha$. Let 
$$
R=\{\pm \epsilon_i\pm \epsilon_j\,|\, 1\leq i\leq j\leq \ell\}\backslash\{0\}
$$ be the corresponding root system realized in $\mathbb R\sp\ell$ with
the canonical basis $\epsilon_1,\dots,\epsilon_\ell$.
Fix simple roots
$$\alpha_1=\epsilon_1-\epsilon_2,\quad\dots,\quad\alpha_{\ell-1}=\epsilon_{\ell-1}-\epsilon_\ell,\quad\alpha_\ell=2\epsilon_\ell$$ 
and let   $ \g=\n_-+ \h + \n_+$ be the corresponding triangular decomposition. Let
$\theta=2\alpha_1+\dots+2\alpha_{\ell-1}+\alpha_\ell=2\epsilon_1$ be the maximal root and
$$
\omega_r= \epsilon_1+\dots+ \epsilon_r, \qquad r=1,\dots,\ell
$$ 
fundamental weights (cf. \cite{H}). 
Fix root vectors $x_\alpha\in\g_\alpha$. % and denote by $\alpha\sp\vee\in\mathfrak h$ dual roots. 
We identify $\mathfrak h$ and $\mathfrak h\sp*$ via the Killing form $\langle\,,\,\rangle$ normalized in 
such a way that  $\langle\theta,\theta\rangle=2$.

\bigskip

%%%%%%%%%%%%%%%%%%%%%%%%%%%%%%%%%%%%%%%%%%%%%%%%%%% Afina Liejeva algebra

Let $\tilde{\mathfrak g}$ be the affine Lie algebra of type $C_\ell\sp{(1)}$ associated  to $\g$,
$$
\gt = \mathfrak{g} \otimes
\mathbb{C}[t,t^{-1}] + \mathbb{C}c, + \mathbb{C} d,
$$
with the canonical central element $c$ and the degree element $d$ (cf. \cite{K}). Let  
$$
\gt=\tilde{\n}_-+ \tilde{\h} + \tilde{\n}_+,
$$
be a triangular decomposition of $\gt$, where $
\tilde{\n}_-=\n_-+ \g\otimes t^{-1}\C [t^{-1}]$, $\tilde{\h}=\h  + \mathbb{C}c + \mathbb{C} d$, $\tilde{\n}_+=\n_++ \g\otimes t\C [t]$.
Denote by $\Lambda_0,\dots,\Lambda_\ell$ fundamental weights of $\gt$.
 
For $x\in{\mathfrak g}$ and $n\in\mathbb Z$ denote by $x(n)=x\otimes t^{n}$ and $x(z)=\sum_{n\in\mathbb Z}
x(n) z^{-n-1}$, where $z$ is a  formal variable. 

Let $L(\Lambda)$ be a standard $\tilde{\mathfrak g}$-module with the highest weight
$$
\Lambda=k_0 \Lambda_0+k_1 \Lambda_1+\dots+k_\ell \Lambda_\ell,
$$
$k_i\in\Z_+$ for $i=0,\dots,\ell$, and fix a highest weight vector $v_\Lambda$.
Denote by
$k=\Lambda(c)$ the level of  $\tilde{\mathfrak g}$-module $L(\Lambda)$,
$
k=k_0 +k_1 +\dots+k_\ell$.

\bigskip

%%%%%%%%%%%%%%%%%%%%%%%%%%%%%%%%%%%%%%%%%%%%%%%%%%% Z-gradacija

Fix the minuscule weight 
$\omega=\omega_\ell=\epsilon_1+\dots+\epsilon_\ell \in\h^*$; then $\langle\omega,\alpha\rangle\in\{-1,0,1\}$ for all $\alpha\in R$ and
define {\it the set of colors} 
$$\Gamma  =
\{\,\alpha \in R \mid \langle\omega,\alpha\rangle = 1\}  =  \{\epsilon_i+\epsilon_j\,|\, 1\leq i \leq j\leq \ell\}.
$$
Write
$$(ij)=\epsilon_i+\epsilon_j\in\Gamma\quad\textrm{and}\quad x_{ij}=x_{\epsilon_i+\epsilon_j}.$$

This gives a $\mathbb Z$-gradation of $\gt$; let ${\mathfrak g}_0 = {\mathfrak h} +
\sum_{\langle\omega,\alpha\rangle=0}\, {\mathfrak g}_\alpha$, then
$$
\gt=\gt_{- 1}+\gt_{0}+\gt_{ 1},
$$
where
$$
\gt_0 = \g_0 \otimes\C [t,t^{-1}] \oplus \C c \oplus \C d,\qquad
\gt_{\pm 1} = \sum_{\alpha \in \pm\Gamma}\, {\mathfrak g}_\alpha \otimes\C [t,t^{-1}].
$$
The subalgebra $\gt_1$ is commutative, and $\g_0$ acts on $\gt_1$ by adjoint action.

\bigskip

%%%%%%%%%%%%%%%%%%%%%%%%%%%%%%%%%%%%%%%%%%%%%%%%%%%%%%%%%%%%%%%%%%%%%% Feigin-Stoyanovsky subspace

{\it Feigin-Stoyanovsky's type subspace} of $L(\Lambda)$ is a $\gt_1$-submodule of $L(\Lambda)$ generated by the highest-weight vector $v_\Lambda$, 
$$
W(\Lambda) = U(\gt_1)\cdot v_\Lambda = U(\gt_1^-)\cdot v_\Lambda\subset L(\Lambda),
$$
where $\gt_1^-=\gt_1 \cap \tilde{\n}_-$.

\bigskip

%%%%%%%%%%%%%%%%%%%%%%%%%%%%%%%%%%%%%%%%%%%%%%%%%%%%%%%%%%%% DC i IC

A monomial $\m\in U(\gt_1^-)= S(\gt_1^-)$,
$$\m=\cdots x_{i'
	j'}(-n-1)^{b_{i'j'}}\cdots x_{i j}(-n)^{a_{ij}}  \cdots,$$% \quad\in U(\gt_1^-)$$ 
is said to satisfy {\em 	difference conditions} for $W(\Lambda)$, $DC$ for short,  if 
for any $n\in\N$ and $i_1 \leq \dots \leq i_t \leq j_t\leq\ldots\leq j_1 \leq i_{t+1} \leq \dots \leq i_{s} \leq j_{s} \leq \ldots\leq j_{t+1}$,
$$ %\label{DC_def}
b_{i_1 j_1}+\ldots+b_{i_t j_t}+a_{i_{t+1} j_{t+1}}+\ldots+a_{i_s
	j_s}\leq k,
$$
where $a_{ij}$'s and $b_{ij}$'s denote exponents of $x_{ij}(-n)$ and
$x_{ij}(-n-1)$ in $\m$, respectively.

A monomial $\m$ satisfies {\em initial conditions} for $W(\Lambda)$, $IC$ for short, if 
for every $i_1 \leq \dots \leq i_t \leq j_t\leq\ldots\leq j_1$,
$$ %\label{IC_def}
a_{i_1 j_1}+\ldots+a_{i_t 	j_t}\leq k_0+k_1+\ldots + k_{j_1-1}
$$
where $a_{ij}$'s denote exponents of $x_{ij}(-1)$ in $\m$.

\begin{tm}[{\cite{BPT}}]%\label{FSbaza_tm}
	The set
	$$\{\m v_\Lambda \,|\, \m \ \textrm{satisfies}\ DC\ \textrm{and}\ IC\ \textrm{for}\ W(\Lambda)\}$$
	is a basis of $W(\Lambda)$.
\end{tm}

\section{Presentation of Feigin-Stoyanovsky's type subspaces}

%%%%%%%%%%%%%%%%%%%%%%%%%%%%%%%%%%%%%%%%%%%%%%%%%%%%%%%%%%%%%%%%%% Relacije DC

Difference conditions are consequences of the adjoint action of $\g_0$ on the vertex-operator relation
$$x_\theta(z)^{k+1}=0,$$ 
or, equivalently, on a family of relations
$$\sum_{\substack{n_1,\dots,n_{k+1}\geq 1 \\ n_1+\dots+n_{k+1}=N}}x_{11}(-n_1)\cdots x_{11}(-n_{k+1})=0,\quad \textrm{for}\ N\geq k+1$$
on $L(\Lambda)$ (cf. \cite{BPT}; see also  \cite{LP}, \cite{MP}, \cite{LL})). 

Root vectors of $\g$ can be chosen so that the action of $\g_0$ on $\g_1$ is given by 
$$
 [x_{-\alpha_i},x_{ij}]=x_{i+1,j}, \quad [x_{-\alpha_j},x_{ij}]=x_{i,j+1},\quad [x_{-\alpha_i},x_{ii}]= 2x_{i,i+1}
$$
(cf. \cite{H}). Then one easily sees that the adjoint action gives the following family of relations on $L(\Lambda)$:
$$
%\label{DC_rel}
\sum_{\substack{n_1+\dots+n_{k+1}=N\\ \{i_1,\dots,i_{k+1},j_1,\dots,j_{k+1} \}=\{1^{m_1},\dots,\ell^{m_\ell}\} }} 
\hspace{-3ex} C_{\bf ij} x_{i_1 j_1}(-n_1)\cdots x_{i_{k+1} j_{k+1}}(-n_{k+1})=0,
$$
for some nonnegative integers $C_{\bf ij}$, where the sum runs over all such partitions ${\bf i},{\bf j}$ of a multiset $\{1^{m_1},\dots,\ell^{m_\ell}\}$, $m_1+\dots+m_\ell=2(k+1)$ (cf. \cite{BPT}).
 
One obtains the difference conditions by finding minimal monomials of these relations,
the so called {\em leading terms} of relations, whose multiples can be excluded from the spanning set 
For this, a linear order on monomials is introduced.
Define a linear order on the set of colors $\Gamma$:
$(i'j')<(ij)$ if $i'>i$ or $ i'=i,\, j'>j$. 
On the {\em set of variables} $\Gamt^-=\{x_\gamma(n) \,|\,\gamma\in\Gamma, n\in \Z_- \}$  define a linear order by
$x_\alpha(n)<x_\beta(n')$ if $n<n'$ or
$n=n',\, \alpha<\beta$. 
For monomials, assume that factors descend from right
to left, then use a lexicographic order (compare factors the greatest to the lowest one).
Order $<$ is compatible with multiplication (see \cite{P1}, \cite{T1}): if $\m_1<\m_2$ then $\m \m_1<\m \m_2$, for  $\m,\m_1, \m_2 \in U(\gt_1^-)$.

\bigskip

%%%%%%%%%%%%%%%%%%%%%%%%%%%%%%%%%%%%%%%%%%%%%%%%%%%%%%%%%%%%%%%%%% Relacije IC

For initial conditions consider decompositions 
%$$\Lambda=\Lambda^{(r)}+\Lambda_{(r)}, \qquad \Lambda^{(r)}=\sum_{t< r } k_t \Lambda_t,\quad \Lambda_{(r)}=\sum_{t\geq r}  k_t \Lambda_t.$$
$$\Lambda=\Lambda^{(r)}+\Lambda_{(r)}, \qquad \Lambda^{(r)}=k_0 \Lambda_0+\dots+k_{r-1}\Lambda_{r-1},\quad \Lambda_{(r)}=k_r \Lambda_r+\dots + k_\ell \Lambda_\ell,$$
for $1\leq r \leq \ell$. 
By $v^{(r)}$ and $v_{(r)}$ denote highest weight vectors of the 
associated standard modules $L(\Lambda^{(r)})$ and $L(\Lambda_{(r)})$ of level
$k^{(r)}=k_0+\dots+k_{r-1} $ and $k_{(r)}=k_{r}+\dots+k_\ell$, respectively. Then $L(\Lambda)$ can be embedded in a tensor product 
$L(\Lambda)\subset L(\Lambda^{(r)}) \otimes L(\Lambda_{(r)})$.  
Since $x_{ij}(-1) v_{\Lambda_r} = 0$ if and only if $j\leq r$ (cf. \cite{BPT}), we have
$$\m (v^{(r)}\otimes v_{(r)}) = (\m v^{(r)}) \otimes v_{(r)}$$
for $\m=x_{i_1 j_1}(-1) \cdots x_{i_t j_t}(-1)$ such that $j_s \leq r$, $1\leq s \leq t$.
Hence, relations between such monomials corresponding to difference conditions for $W(\Lambda^{(r)})$ automatically become 
relations in $L(\Lambda)$ (cf. \cite{BPT}).
This gives the following family of relations on $L(\Lambda)$: 
\begin{equation}
\label{IC_rel}
\sum_{\{i_1,\dots,i_{k+1},j_1,\dots,j_{k+1} \}=\{1^{m_1},\dots,r^{m_r}\} } 
\hspace{-3ex} C_{\bf ij} x_{i_1 j_1}(-1)x_{i_2 j_2}(-1)\cdots x_{i_{k^{(r)}+1} j_{k^{(r)}+1}}(-1)=0,
\end{equation}
for some nonnegative integers $C_{\bf ij}$, where the sum runs over all such partitions ${\bf i}, {\bf j}$ of a multiset $\{1^{m_1},\dots,r^{m_r}\}$, $m_1+\dots+m_r=2(k^{(r)}+1)$.

Alternatively, for $r\geq 2$ let $\g_{(r)}\subset \g_0$ be the subalgebra generated by elements $x_{\pm\alpha_t}$, $1\leq t < r$.
Start from a relation 
$$x_{11}(-1)^{k^{(r)}+1} (v^{(r)}\otimes v_{(r)}) = (x_{11}(-1)^{k^{(r)}+1} v^{(r)})\otimes v_{(r)} = 0.$$
Now the adjoint action of  $\g_{(r)}$ on the above relation gives  relations \eqref{IC_rel}. For $r=1$,  relations \eqref{IC_rel} come down to only one relation
$$x_{11}(-1)^{(k_0+1)} (v^{(1)}\otimes v_{(1)})=0.$$

%%%%%%%%%%%%%%%%%%%%%%%%%%%%%%%%%%%%%%%%%%%%%%%%%%%%%%%%%%%% Prezentacija

\bigskip

Recall that Feigin-Stoyanovsky's type subspace $W(\Lambda)$ is 
$$W(\Lambda)=U(\gt_1^-)\cdot v_\Lambda.$$
Since $\gt_1$ is commutative, universal enveloping algebra of $\gt_1^-$ is isomorphic to a polynomial algebra $\C[\Gamt^-]$. Hence, there is a 
surjection $$f_\Lambda:\C[\Gamt^-]\to W(\Lambda), \qquad f:\m \to \m\cdot v_\Lambda.$$
We want to describe the kernel of this map, $\ker f_\Lambda\subset \C[\Gamt^-]$, so that
$$W(\Lambda) \simeq \C[\Gamt^-]/\ker f_\Lambda, $$
as vector spaces.

\begin{tm}
Let $J_\Lambda\subset\C[\Gamt^-]$ be the ideal generated by the following sets
\begin{eqnarray*} 
%\label{Pres_eq1} 
U(\g_0)\cdot\left(\sum_{\substack{n_1,\dots,n_{k+1}\geq 1 \\ n_1+\dots+n_{k+1}=N}}x_{11}(-n_1)\cdots x_{11}(-n_{k+1})\right),& &  \textrm{for}\ N\geq k+1,\\
%\label{Pres_eq2} 
U(\g_{(r)})\cdot x_{11}(-1)^{k^{(r)}+1},& &  \textrm{for}\ r=2,\dots,\ell,\\
%\label{Pres_eq3} 
x_{11}(-1)^{k_0+1}.& &   
\end{eqnarray*}
Then $\ker f_\Lambda = J_\Lambda$.
\end{tm}

\begin{dokaz}
The generators of $J_\Lambda$ are exactly the elements that appear on the left side of relations on $W(\Lambda)$. Therefore they
lie in the kernel of $f_\Lambda$, and we can factorize $f_\Lambda$
to a quotient map
$$	\bar{f_\Lambda}: \C[\Gamt^-]/J_\Lambda \to  W(\Lambda).	$$ 
This map is clearly a surjection, since $f_\Lambda$ is a surjection. 

We can imitate the proof for the spanning set for $W(\Lambda)$ (cf. Proposition 2 and 4 in \cite{BPT}) to reduce the spanning set for $\C[\Gamt^-]/J_\Lambda$. 
Consider the generators of $J_\Lambda$ and identify the minimal monomial inside each one; their multiples can be excluded from the spanning set. 
Like in \cite{BPT},  we get 
	$$\mathcal{B}=\{\m\,|\,\m\ \textrm{satisfies DC and IC for}\ W(\Lambda)\}$$
as a spanning set of $\C[\Gamt^-]/J_\Lambda$.
	
To see that $\bar{f_\Lambda}$ is an injection, note that $\bar{f_\Lambda}$ maps $\mathcal{B}$ bijectively onto 
$$\{\m v_\Lambda\,|\,\m \ \textrm{satisfies DC and IC for}\ W(\Lambda)\}\subset W(\Lambda),$$
which is a basis of $W(\Lambda)$.  This means that $\mathcal{B}$ is also linearly
independent. Hence $\bar{f_\Lambda}$ maps a basis of $\C[\Gamt^-]/J_\Lambda$ onto a basis 
of $W(\Lambda)$ and therefore $\bar{f_\Lambda}$ is a bijection.
\end{dokaz}

\end{document}